\documentclass[a4paper, twoside,11pt]{article} %Классовые опции titlepage, draft, fleqn, leqno, twoside,twocolumn

\usepackage[utf8]{inputenc}
\usepackage[english, russian]{babel}
\usepackage{amsmath,amsfonts,amssymb,amscd,euscript,amsthm}
\usepackage{ upgreek } %% for \count -- first infinite ordinal

\theoremstyle{definition}

\theoremstyle{plain}
\newtheorem{teo}{Theorem}[section]
\newtheorem{cor}{Corollary}[section]
\newtheorem{lem}{Lemma}[section]
\theoremstyle{remark}
\newtheorem{rem}{Remark}[subsection]

\title{Unlocking of predicate: application to non-anticipating selections}  % Declares the document's title.

\author{D.\,A.~Serkov}      % Declares the author's name.

\date{Krasovskii Institute of Mathematics and Mechanics, UB RAS;\\
Ural Federal University named after B.N. Yeltsin;\\
E-mail: serkov@imm.uran.ru}

\newcommand{\beq}{\begin{equation}}
\newcommand{\beqnt}{\begin{equation}\notag}
\newcommand{\eeq}{\end{equation}}
\newcommand{\bpr}{\begin{proof}}
\newcommand{\epr}{\end{proof}}
\newcommand{\fref}[1]{{\rm(\ref{#1})}}                      %%%%%%%%%% formula ref      %%%%%
\newcommand{\mydef}{\ensuremath{\stackrel{\text{\mdef}}{=}}}% "равно по определению" в тексте
\newcommand{\mdef}{\textup{def}}% "равно по определению" в формуле
\newcommand{\mdefeq}{\stackrel{\mdef}{\myeqv}}% "эквив. по определению" в формуле
\renewcommand{\le}{\leqslant}% Меньше равно
\newcommand{\myle}{\mathrel\preccurlyeq}% знак частичного порядка
\newcommand{\myLe}{\mathrel\sqsubseteq}% 2-ой знак ч.п.
\newcommand{\myor}{\ensuremath{\lor}}% или
\newcommand{\myand}{\ensuremath{\&}}% и
\newcommand{\myno}{\ensuremath{\neg}}% не
\newcommand{\myimp}{\ensuremath{\Rightarrow}}% имп
\newcommand{\myeqv}{\ensuremath{\Leftrightarrow}}% экв
\newcommand{\myemp}{\ensuremath{\varnothing}}% пустое мн.
\newcommand{\myll}{\ensuremath{\forall}}% квантор всеобщности
\newcommand{\myxst}{\ensuremath{\exists}}% квантор существования
\newcommand{\icP}{\ensuremath{\EuScript P}}% булеан
\newcommand{\res}[2]{\ensuremath{(#1\,|\,#2)}}% сужение ф-ции #1 на #2
\newcommand{\sres}[2]{\ensuremath{\text{\large\textbf{\textup(}}#1\,|\,#2\text{\large\textbf{\textup)}}}}% сужение множества #1 на #2
\newcommand{\ares}[2]{\ensuremath{[#1\,|\,#2]}}% сужение м/з-автом #1 на #2
\newcommand{\fix}[1]{\ensuremath{{\mathbf{Fix}({#1})}}}%%%%%   множество неподвижных точек отображения
\newcommand{\Ntst}{\ensuremath{\mathcal X}}%%% тестовое семейство неупреждаемости
\newcommand{\spri}[3]{\ensuremath{#1\res{#2}{#3}}}% росток в #1 функции #2 со множества #3
\newcommand{\hspri}[3]{\ensuremath{[#1]\res{#2}{#3}}}% общий след откликов м\ф #1 на росток в омега помехи #2 со множества #3
\newcommand{\puI}{\ensuremath{\mathcal I}}%%%%%%% множество индексов в произведении
\newcommand{\puJ}{\ensuremath{\mathcal J}}%%%%%%% множество индексов в  конъюнкции
\newcommand{\UM}[1]{\ensuremath{{\mathfrak{UM}({#1})}}}%%%%%   множество разм. отображений пред.#1
\newcommand{\PRE}[1]{\ensuremath{{\mathfrak{PR}({#1})}}}%%%%%   множество пред. на #1
\newcommand{\PTRUE}{\ensuremath{{\mathfrak{T}}}}%%%%%   TRUE
\newcommand{\PFALS}{\ensuremath{{\mathfrak{F}}}}%%%%%   FALSE
\newcommand{\icI}{\ensuremath{I}}% время
\newcommand{\icX}{\ensuremath{X}}% прос-во
\newcommand{\icOm}{\ensuremath{\Omega}}% неопред.
\newcommand{\icD}{\ensuremath{D}}% D=IxX
\newcommand{\icC}{\ensuremath{\mathbf C}}% траектор.
\newcommand{\icY}{\ensuremath{Y}}% прос-во неопред.
\newcommand{\icM}{\ensuremath{\mathbf M}}% м/ф (\icOm,\icC)
\newcommand\muf{\ensuremath{\mathcal{M}}} %%%%%%%%%% заданная мультифункция
\newcommand{\icORD}{\ensuremath{\mathbf{ORD}}}% семейство ординалов
\newcommand{\LE}{\mathbf{LE}}%%%%%   нижняя половинка TRUE
\newcommand{\orle}{\ensuremath{\mathrel\preccurlyeq}}% порядок на ординалах
\newcommand{\cttv}{restrictive}% сужающее
\newcommand{\carle}{\ensuremath{\mathrel\text{\guilsinglleft=}}}% порядок на кардиналах
\newcommand{\carles}{\ensuremath{\mathrel\text{\guilsinglleft}}}% строгий порядок на кардиналах
\newcommand{\CARD}[1]{\ensuremath{|#1|}}% множесво множеств равномощных #1 
\newcommand{\CRDP}[1]{\ensuremath{|#1|^+}}% наим кард более мощ чем #1 
\newcommand{\fiord}{\ensuremath{\upomega}}% первый бесконечный предельный ординал

\begin{document}
%\rissue

\maketitle

%%%%%%%%%%%%%%%%%%%%%%%%%%%%%%%%%%%%%%%%%%%%%%%%%%%%%%%%%%
\section{Introduction}

We consider an approach to constructing a non-anticipating selection of a multivalued mapping; such a problem arises in control theory under conditions of uncertainty.
The approach is called ``unlocking of predicate'' and consists in the reduction of finding the truth set of a predicate to searching fixed points of some mappings.
Unlocking of predicate gives an extra opportunity to analyze the truth set and to build its elements with desired properties.

This concept is used in many fields of mathematics:
in differential equations and differential inclusions;
in game theory, when studying the saddle points (see \cite{Kakutani1941}) and the Nash equilibria (see \cite{Nash-Annal-1951, Nikaido-1954});
in dynamic games, when constructing the stable (weakly invariant) sets (see \cite{Chentsov75DANet, Chentsov76DAN226et}) and non-anticipating selections of multivalued mappings (see \cite{Chentsov98DUPUe, Chentsov99DUPUe}).
However, in all the above cases ``unlocking mappings'' are presented as a ready--made product:
a method for constructing an ``unlocking mapping'' has remained beyond the consideration.

In this article, we outline how to build ``unlocking mappings'' for some general types of predicates:
we give a formal definition of the predicate unlocking operation, the rules for the construction and calculation of ``unlocking mappings'' and their basic properties.
As an illustration, we routinely construct two unlocking mappings for the predicate ``be non-anticipating mapping'' and then on this base we provide the expression for the greatest non-anticipating selection of a given multifunction.
This work continues \cite{Ser_UDSU2016e} where the procedure for the predicate "be Nash equilibrium" is presented.
%The described approach is far from universality, but at least can be used for the above mentioned positive examples.

%%%%%%%%%%%%%%%%%%%%%%%%%%%%%%%%%%%%%%%%%%%%%%%%%%%%
\section{Notation and definitions}

1. 
Hereinafter, we use the set--theoretic symbols (quantifiers, propositional bundles, \myemp\ for the empty set); 
$\mydef$\ for the equality by definition; 
%"\mdef"\ replaces the phrase "by definition"; 
$\mdefeq$ for the equivalence by definition.
We accept the axiom of choice. 
A set consisting of sets is called a family. 
By $\icP(T)$ (by $\icP'(T)$), we denote the family of all (all nonempty) subsets of an arbitrary set $T$; the family $\icP(T)$ is also called Boolean of the set $T$.
If $A$ and $B$ are non--empty sets, then $B^A$ is the set of all functions from the set $A$ to the set $B$ (see \cite{KurMos1970et}). 
If $f\in B^A$ and $C\in\icP'(A)$, then $\res{f}{C}\in B^C$ is the restriction of $f$ to the set $C$: $\res{f}{C}(x)\mydef f(x)$ $\myll x\in C$. 
We denote the image of the set $C\in\icP(A)$ under the function $f$ by $f(C)$: $f(C)\mydef\{f(x): x\in C\}$.
When $f\in\icP(B)^A$, $f$ is called a multivalued function or multifunction (m/f) from $A$ in $B$.
The term ``mapping'' means a function or m/f.
In case $F\in\icP'(B^A)$, we denote $\sres{F}{C}\mydef\{\res{f}{C}:f\in F\}$.
If $f\in B^A$, we denote by $f^{-1}$ the m/f from $B$ into $A$ defined by the rule
$$
f^{-1}(b)\mydef
\begin{cases}
\{a\in A\mid b=f(a)\},& b\in f(A),\\
\myemp,& b\not\in f(A)
\end{cases}
\qquad \myll b\in B.
$$
We call the m/f $f^{-1}$ \emph{inverse mapping}.
If $f\in\icP(B)^A$, i.e. $f$ is a m/f, we define $f^{-1}$ by
$$
f^{-1}(b)\mydef
\begin{cases}
\{a\in A\mid b\in f(a)\},& b\in\bigcup f(A),\\
\myemp,& b\not\in\bigcup f(A),
\end{cases}
\qquad \myll b\in B.
$$
For any $f\in X^X$ we denote by $\fix{f}$ the set of all fixed points of $f$: $\fix{f}\mydef\{x\in X\mid f(x)=x\}$.
In the case when $f$ is a m/f, the set $\fix{f}$ is defined by: $\fix{f}\mydef\{x\in X\mid x\in f(x)\}$.

2. 
{\em A predicate} $P$ on a non-empty set $X$ is identified with the same name function from $\{0,1\}^X$.
We say that $x\in X$ satisfies the predicate $P$ and write it down by $P(x)$ iff $P(x)=1$.
The set of all $x\in X$ satisfying the predicate $P$ is called \emph{the set of truth} (of the predicate $P$).
Following the definition of an inverse mapping, we denote this set by $P^{-1}(1)$.
The set of all predicates on $X$ is denoted by $\PRE{X}$.
We denote by $\PTRUE$ ($\PFALS$) the predicate on $X$ defined by: $\PTRUE^{-1}(1)=X$ ($\PFALS^{-1}(0)=X$).
Hence, for any $P\in\PRE{X}$, the equalities $P=\PTRUE\myand P=\PFALS\myor P$, where ``\myand'' (``\myor'') denotes logical ``and'' (``or''),  are valid.

We call \emph{unlocking of predicate} $P$ the operation of constructing a mapping $\mathcal F_P\in\icP(X)^X\cup X^X$  that satisfies the condition
\beq\label{P-unlock}
\fix{\mathcal F_P}= P^{-1}(1).
\eeq
The mapping $\mathcal F_P$ with property \fref{P-unlock}, is called \emph{unlocking mapping} (for the predicate $P$).
Denote by \UM{P} the set of all unlocking m/f for the predicate $P$.
Thus, $\UM{P}\in\icP'(\icP(X)^X)$.
The formal exclusion of functions (the set $X^X$) from \UM{P}\ is dummy, because every function $f$ satisfying $\fix{f}=P^{-1}(1)$ is represented by the m/f $F_f$ in \UM{P}: $F_f(x)\mydef\{f(x)\}$ $\myll x\in X$.
So, for a function $f$ we write down $f\in\UM{P}$ keeping in mind the inclusion $F_f\in\UM{P}$.

3. 
For any set $X\neq\myemp$ and a partial ordering relation $\myle\in\icP'(X\times X)$, we denote by $(X, \myle)$ the corresponding partially ordered set (poset). 
A set $C\subset X$ is called a \emph{chain} if it is totally ordered by $\myle$: $(x\myle y)\myor(y\myle x)$ $\myll x,y\in C$.
In particular, \myemp\ is a chain. 
Following \cite{Markow1976}, we call a poset $(X,\myle)$ a\emph{chain--complete} poset if there exists the greatest lower bound $\inf C\in X$ for any chain $C\subset X$. 
In particular, every chain--complete poset $(X, \myle)$ has the greatest element $\top\in X$ (the greatest lower bound of the empty chain), and, thus, it is not empty. 
For $Y\in\icP(X)$, we denote by $\top_Y$ and $\bot_Y$ the greatest and the least elements of the set $Y$, respectively, if they exist. 
A poset is called a \emph{complete lattice} iff any subset has the greatest and the least elements.
So, any complete lattice is a chain--complete poset.
Let $(X, \myle)$ be a non-empty poset and $f\in X^X$. 
The function $f$ is called \emph{\cttv} if $f(x)\myle x$ for every $x\in X$ . 
The function $f$ is called \emph{isotone} if the implication $(x\myle y)\myimp(f(x)\myle f(y))$ holds for all $x,y\in X$. 

4.
Denote the class of ordinals by \icORD. 
For a set $X$, we denote by $\CARD{X}$ the least ordinal that is equipotent to the set $X$ (the cardinal of $X$).
The relation of order (strict order) on the class of cardinals is denoted by \carle\ (\carles).
For any set $H$, let $\CRDP{H}\in\icORD$ be the least ordinal among the ordinals $\eta$ with the property $\CARD{H}\carles\CARD{\eta}$.

%%%%%%%%%%%%%%%%%%%%%%%%%%%%%%%%%%%%%%%%%%%%%%%%%%%%%
\section{Calculus of unlocking mappings}\label{sec-P-unlock}

\subsection{The order, restrictions and logical operations}

1. 
Let $X$ be a nonempty set.
On the set $\icP(X)^X$, we introduce the partial order $\myle$, assuming that 
$
(g\myle f)\mdefeq(g(x)\subset f(x)\ \myll x\in X)\qquad\myll f,g\in\icP(X)^X
$
Then we have the equivalence
$%\beq\label{le-rev}
(g\myle f)\myeqv(g^{-1}\myle f^{-1}).
$%\eeq
The poset $(\icP(X)^X,\myle)$ is a complete lattice.
It is also easy to check that, for any $P\in\PRE{X}$, the poset $(\UM{P},\myle)$ forms a complete sublattice (a subset being a complete lattice) in $(\icP(X)^X,\myle)$ and the equalities are true: 
$$
\top_{\UM{P}}(x)=
\begin{cases}
X,&P(x),\\
X\setminus\{x\},&\myno P(x),
\end{cases}
\qquad
\bot_{\UM{P}}(x)=
\begin{cases}
\{x\},&P(x),\\
\myemp,&\myno P(x).
\end{cases}
$$
In particular, for the predicates $\PTRUE$, $\PFALS$, the relations $\top_\PTRUE(x)=X$, $\bot_\PTRUE(x)=\{x\}$, $\top_\PFALS(x)=X\setminus\{x\}$, and $\bot_\PFALS(x)=\myemp$ are valid for all $x\in X$.
By definition, we have $\top_{\UM{P}}=\top_{\UM{P}}^{-1}$, $\bot_{\UM{P}}=\bot_{\UM{P}}^{-1}$.

%%%%%%%%%%%%%%%%%%%%%%%%%%%%%%%%%%%%%%%%%%%%%%%%%%%%%%%
\begin{lem}\label{lem-UMP}
For all $f\in\icP(X)^X$, the relations $(f\myle\top_{\UM{P}})\myimp(\fix{f}\subset P^{-1}(1))$, $(\bot_{\UM{P}}\myle f)\myimp(P^{-1}(1)\subset\fix{f})$ are fulfilled.
Consequently,
\beqnt%\label{UMP}
\UM{P}=\{f\in\icP(X)^X\mid\bot_{\UM{P}}\myle f\myle\top_{\UM{P}}\},
\qquad
(f\in\UM{P})\myeqv (f^{-1}\in\UM{P}).
\eeq
\end{lem}

2. 
For any $\phi\in\icP(X)^X$ and $Y\in\icP'(X)$, we denote by $\ares{\phi}{Y}$ the following mapping $\ares{\phi}{Y}(y)\mydef Y\cap\phi(y)$ $\myll y\in Y$.
Recall that the restriction $\res{P}{Y}\in\PRE{Y}\mydef\{0,1\}^Y$ of $P\in\PRE{X}$ is defined by $\res{P}{Y}(y)\mydef P(y)$, $\myll y\in Y$.
%%%%%%%%%%%%%%%%%%%%%%%%%%%%%%%%%%%%%%%%%%%%%%%%%%%%%%%
\begin{lem}\label{lem-P-res}
For all $P\in\PRE{X}$, $Y\in\icP'(X)$ the equalities $\UM{\res{P}{Y}}=\{\ares{\phi}{Y}:\phi\in\UM{P}\}$ are valid.
\end{lem}
%%%%%%%%%%%%%%%%%%%%%%%%%%%%%%%%%%%%%%%%%%%%%%%%%%%%%%%

3. 
The following lemma provides unlocking mappings for some expressions of propositional logic.

%%%%%%%%%%%%%%%%%%%%%%%%%%%%%%%%%%%%%%%%%%%%%%%%%%%%%%%
\begin{lem}\label{lem-PQ-unlock}
If $P,Q\in\PRE{X}$, then the equalities are valid:
$$%\begin{multline*}
\UM{\myno P}=\{f\in\icP(X)^X\mid\myxst g\in\UM{P}: f(x)=X\setminus g(x)\ \myll x\in X\},%\label{notP-unlock}
$$%\end{multline*}
$$%\begin{multline*}
\UM{P\myand Q}=\{f\in\icP(X)^X\mid\myxst g\in\UM{P}\myxst q\in\UM{Q}: f(x)=g(x)\cap q(x)\ \myll x\in X\},%\label{PandQ-unlock}
$$%\end{multline*}
$$%\begin{multline*}
\UM{P\myor Q}=\{f\in\icP(X)^X\mid\myxst g\in\UM{P}\myxst q\in\UM{Q}: f(x)=g(x)\cup q(x)\ \myll x\in X\}.%\label{PorQ-unlock}
$$%\end{multline*}
\end{lem}
%%%%%%%%%%%%%%%%%%%%%%%%%%%%%%%%%%%%%%%%%%%%%%%%%%%%%%%

Using the above relations, one can construct unlocking mappings for a variety of other propositional logic expressions.

%%%%%%%%%%%%%%%%%%%%%%%%%%%%%%%%%%%%%%%%%%%%%%%%%%%%%%%
\begin{cor}\label{lem-FTP}
For all  $P\in\PRE{X}$, $f\in\UM{P}$, $T\in\UM{\PTRUE}$, and $F\in\UM{\PFALS}$, the m/f $f_T, f_F\in\icP(X)^X$ defined by $f_T(x)\mydef T(x)\cap f(x)$, $f_F(x)\mydef F(x)\cup f(x)$ $\myll x\in X$,
are unlocking m/f for the predicate $P$: $f_T, f_F\in\UM{P}$.
In addition, we have the relations:
\begin{gather*} 
\bot_{\UM{P}}(x)=\bot_{\UM{\PTRUE}}(x)\cap f(x)=\{x\}\cap f(x),\label{UM-P-bot}\\
\top_{\UM{P}}(x)=\top_{\UM{\PFALS}}(x)\cup f(x)=(X\setminus\{x\})\cup f(x).\label{UM-P-top}
\end{gather*}
\end{cor}
%%%%%%%%%%%%%%%%%%%%%%%%%%%%%%%%%%%%%%%%%%%%%%%%%%%%%

4.
Lemma \ref{lem-f-contr} is based on the corollary \ref{lem-FTP} and allows us  to construct an unlocking function from a given unlocking m/f in the case when $X$ is an ordered set.
Let $(X,\le)$ be a poset and the m/f $\LE_X\in\icP(X)^X$ is defined by 
\beq\label{lex}
\LE_X(x)\mydef\{y\in X\mid y\le x\}.
\eeq
Notice that $\LE_X\in\UM{\PTRUE}$.
%%%%%%%%%%%%%%%%%%%%%%%%%%%%%%%%%%%%%%%%%%%%%%%%%%%%%%%
\begin{lem}\label{lem-f-contr}
Let $(X,\myle)$ be a nonempty poset, $P\in\PRE{X}$, and $f\in\UM{P}$. 
Let $G\in\icP(X)^X$ be defined by 
\beq\label{def-G}
G(x)\mydef\LE_X(x)\cap f(x)\qquad x\in X,
\eeq
$Y\mydef\{y\in X\mid G(y)\neq\myemp\}$, and the function $g\in Y^Y$ be defined by 
\beqnt%\label{def-g}
g(x)\mydef
\begin{cases}
\top_{G(x)},&\myxst\top_{G(x)},\\
y\in G(x),&\myno\myxst\top_{G(x)},
\end{cases}
\qquad x\in Y.
\eeq
Then $g$ is \cttv\ on $(Y,\myle)$ and $\fix{g}=P^{-1}(1)$.
\end{lem}

%%%%%%%%%%%%%%%%%%%%%%%%%%%%%%%%%%%%%%%%%%%%%%%%%%%%%%%

%%%%%%%%%%%%%%%%%%%%%%%%%%%%%%%%%%%%%%%%%%%%%%%%%%%%%%%%%%%%%%%%%%%%%
\subsection{Unlocking the conjunction of predicates defined on a product}

The conjunction of a set of predicates is an important particular case. 
Using this peculiarity, lemma \ref{lem-FP-unlock} gives the construction of the corresponding unlocking m/f.

Let $\puI$, $(X_\iota)_{\iota\in\puI}$ be non--empty sets and 
\beq\label{setX}
X\mydef\prod_{\iota\in\puI}X_\iota.
\eeq
We call an element $x\in X$ a \emph{tuple from $X$} (or simply a \emph{tuple} if the set $X$ is fixed) and denote the $\iota$--th element of the tuple $x$  by $x_\iota$: $x_\iota\mydef\res{x}{\{\iota\}}\in X_\iota$.
Denote by $(y,x_{-\iota})$ the tuple from $X$ resulted from the tuple $x\in X$ by substituting the element $y\in X_\iota$ at the position of  $x_\iota$:
\beqnt%\label{subst}
(y,x_{-\iota})_\jmath\mydef
\begin{cases}
y,&\jmath=\iota,\\
x_\jmath,&\jmath\in\puI\setminus\{\iota\},
\end{cases}
\qquad\myll x\in X\ \myll y\in X_\iota\ \myll\iota\in\puI.
\eeq

Let a family of predicates $P_\jmath\in\PRE{X}$, $\jmath\in\puJ$ on the product $X$ be given. 
Let the predicate $P\in\PRE{X}$ have the form $P(x)\mdefeq(P_\jmath(x)\ \myll\jmath\in\puJ)$ $x\in X$.
Let $\CARD{\puJ}\carle\CARD{\puI}$ and $q\in\puI^\puJ$ be the corresponding injection of \puJ\ into \puI.
Define the m/f $\mathcal F_P\in\icP(X)^X$ as follows:
\beq\label{Pun}
\mathcal F_P(x)\mydef\prod_{\iota\in\puI}\mathcal B_\iota(x)\qquad\myll x\in X,
\eeq
where m/f $\mathcal B_{\iota},\mathcal B_{\iota\jmath}\in\icP(X_\iota)^X$ are defined by
\beq\label{BPj}
\mathcal B_\iota(x)\mydef
\begin{cases}
\mathcal B_{\iota q^{-1}(\iota)}(x)((y,x_{-\iota}))\},&\iota\in q(\puJ),\\
X_\iota,&\iota\not\in q(\puJ),
\end{cases}
\quad \mathcal B_{\iota\jmath}(x)\mydef\{y\in X_\iota\mid P_\jmath((y,x_{-\iota}))\},\quad x\in X,\ \iota\in\puI,\ \jmath\in\puJ.
\eeq
%%%%%%%%%%%%%%%%%%%%%%%%%%%%%%%%%%%%%%%%%%%%%%%%%%%%%%%
\begin{lem}\label{lem-FP-unlock}
$\mathcal F_P\in\UM{P}$.
\end{lem}
%%%%%%%%%%%%%%%%%%%%%%%%%%%%%%%%%%%%%%%%%%%%%%%%%%%%%%%

%%%%%%%%%%%%%%%%%%%%%%%%%%%%%%%%%%%%%%%%%%%%%%%%%%%%%%
\section{The greatest non-anticipating selection}
\label{Nansel-unlock}

In \cite{Chentsov98DUPUe, Chentsov99DUPUe}, the representation of non-anticipating selections of a m/f as the set of fixpoins of a function (noted by $\Gamma$) is provided.
In other words, the unlocking of predicate ``be non-anticipating selection'' is fulfilled.
At the same time, the process of constructing the function $\Gamma$ remained out of consideration. 
In this section, we carry out the process using constructions from \cite{Chentsov98DUPUe, Chentsov99DUPUe} and relations from the section \ref{sec-P-unlock}.

%%%%%%%%%%%%%%%%%%%%%%%%%%%%%%%%%%%%%%%%%%%%%%%%%%%%
\subsection{Notation and definitions}

Hereinafter, we fix  $\icD\mydef\icI\times\icX$, where \icI and \icX\ are non-empty sets.
Select the set $\icC\in\icP'(\icX^\icI)$ whose elements are considered as ``realizations of control actions''. 
So, the sets \icI\ and \icX\ are analogues of time and state space respectively.
Then we select and fix the sets \icY\ and $\icOm\in\icP'(\icY^\icI)$. 
Elements of $\icOm$ are treated as ``realizations of uncertainty factors''. 
Let $\icM\mydef\icP(\icC)^\icOm$ denote the set of all m/f from \icOm\ into \icC: $\alpha(\omega)\subset\icC$ for any $\omega\in\icOm$, $\alpha\in\icM$.

The partial order $\myLe$ on $\icM$ is defined by
$$
(\phi\myLe\psi)\mdefeq(\phi(\omega)\subset\psi(\omega)\quad\myll\omega\in\icOm)\qquad\myll\phi,\psi\in\icM.
$$
One can verify that the poset $(\icM,\myLe)$ is a complete lattice.
For any $\phi,\psi\in\icM$, we call m/f $\phi$ a \emph{selection} of $\psi$ iff $\phi\myLe\psi$.

Let $\Ntst\in\icP(\icI)$ be a non--empty set. 
A m/f $\phi\in\icM$ is called {\em $\Ntst$--non-anticipating}, iff 
\beq\label{Nt-nant}
(\omega'\in\spri{\icOm}{\omega}{A})\myimp\left(\sres{\phi(\omega)}{A}\subset\sres{\phi(\omega')}{A}\right)\qquad \myll A\in\Ntst,\ \myll\omega,\omega'\in\icOm.
\eeq

\begin{rem}
Due to the equivalence 
$$
(\omega'\in\spri{\icOm}{\omega}{A})\myeqv(\omega\in\spri{\icOm}{\omega'}{A})\myeqv(\res{\omega}{A}=\res{\omega'}{A})\qquad \myll A\in\Ntst,\ \myll\omega,\omega'\in\icOm,
$$
implications \fref{Nt-nant} are equivalent to the relations  
\beqnt%\label{Pna2}
(\res{\omega}{A}=\res{\omega'}{A})\myimp\left(\sres{\phi(\omega)}{A}=\sres{\phi(\omega')}{A}\right)\qquad \myll A\in\Ntst,\ \myll\omega,\omega'\in\icOm,
\eeq
which are usually considered as the definition of non-anticipating property.
\end{rem}

Fix the family $\Ntst$ and a m/f $\muf\in\icM$.
Our aim is to find the greatest in $(\icM,\myLe)$ $\Ntst$-non-anticipating selection of the m/f $\muf$.
So, we should find a m/f $\phi\in\icM$, satisfying condition \fref{Nt-nant}, the inequality $\phi\myLe\muf$, and such that the relation $\beta\myLe\phi$ is valid for any $\beta\in\icM$ satisfying \fref{Nt-nant} and the inequality $\beta\myLe\muf$.

For the analysis of the problem above, we define the predicate $P_{na}\in\PRE{\icM}$ ``be $\Ntst$-non-anticipating mapping'' by
\beq\label{Pna}
P_{na}(\phi)\mdefeq\bigl((\omega'\in\spri{\icOm}{\omega}{A})\myimp\left(\sres{\phi(\omega)}{A}\subset\sres{\phi(\omega')}{A}\right)\ \myll A\in\Ntst\ \myll\omega,\omega'\in\icOm\bigr)\qquad\myll\phi\in\icM,
\eeq
and introduce some new notation. 
For arbitrary $A\in\Ntst$, $\Psi\subset\icOm$, $\omega\in\icOm$, $H\subset\icC$, $h\in\icC$, and $\phi\in\icM$ we set
\beq\label{ospring}
\spri{\Psi}{\omega}{A}\mydef\{\nu\in\Psi\mid\res{\nu}{A}=\res{\omega}{A}\},\qquad\spri{H}{h}{A}\mydef\{f\in H\mid\res{f}{A}=\res{h}{A}\},
\eeq
\beqnt%\label{ospring2}
\spri{\Psi}{-\omega}{A}\mydef\spri{\Psi}{\omega}{A}\setminus\{\omega\}
\eeq
\beq\label{ctrac}
\hspri{\phi}{\omega}{A}\mydef\bigcap_{\nu\in\spri{\icOm}{\omega}{A}}\sres{\phi(\nu)}{A},
\eeq
\beq\label{ctrac2}
\hspri{\phi}{-\omega}{A}\mydef\bigcap_{\nu\in\spri{\icOm}{-\omega}{A}}\sres{\phi(\nu)}{A}.
\eeq

%%%%%%%%%%%%%%%%%%%%%%%%%%%%%%%%%%%%%%%%%%%%%%%%%%%%%%%%
\subsection{Unlocking the predicate ``be \Ntst-non-anticipating mapping''}\label{nonant-unlock}

It follows directly from definition \fref{Pna} that $P_{na}$ is the conjunction of the family $\{P_\omega\mid\omega\in\icOm\}$, where
\beq\label{Pof}
P_\omega(\phi)\mdefeq\bigl((\omega'\in\spri{\icOm}{\omega}{A})\myimp\left(\sres{\phi(\omega)}{A}\subset\sres{\phi(\omega')}{A}\right)\ \myll A\in\Ntst\bigr)\qquad\myll\omega\in\icOm,\ \myll\phi\in\icM.
\eeq
Then we transform \fref{Pof} using notation \fref{ctrac}.
\begin{lem}
\beq\label{Nt-nant-eq-hspring3}
P_\omega(\phi)\myeqv\left(\hspri{\phi}{\omega}{A}=\sres{\phi(\omega)}{A}\quad \myll A\in\Ntst\right)\qquad\myll\omega\in\icOm,\ \myll\phi\in\icM.
\eeq
\end{lem}

So, the predicate $P_{na}$ has the form (see \fref{Pna})
\beqnt%\label{PU}
P_{na}(\phi)\myeqv\left(P_\omega(\phi)\ \myll\omega\in\icOm\right)\myeqv\left(\hspri{\phi}{\omega}{A}=\sres{\phi(\omega)}{A}\quad \myll A\in\Ntst\ \myll\omega\in\icOm\right)\qquad\myll\phi\in\icM.
\eeq

According to scheme \fref{setX} -- \fref{Pun}, we represent $\icM$ as the product of $\icOm$ copies of the set $\icP(\icC)$. 
By the definitions the index set in the conjunction representing $P_{na}$ coincides with the one in the product representing $\icM$.
Hence, the injection $q$ in \fref{BPj} can be chosen as the identity mapping.
Then we have 
$$
\puI\mydef\icOm,\quad X_\iota\mydef X_\omega\mydef\icP(\icC),\quad \omega\in\icOm,
$$
$$
\icM\mydef X\mydef\prod_{\iota\in\puI}X_\iota\mydef\prod_{\omega\in\icOm}\icP(\icC),\quad\mathcal B_\iota\mydef\mathcal B_\omega\in\icP(\icP(\icC))^\icM,
$$
$$
\mathcal F_{P_{na}}(\phi)\mydef\prod_{\omega\in\icOm}\mathcal B_\omega(\phi)\in\icP(\icP(\icC))\mydef\icP(\icM)^\icM.
$$
We provide this list of ``actors and performers'' for the convenience of tracking scheme \fref{setX}--\fref{Pun}.
%In the next calculations we won't switch to the notation from point \ref{sec-P-unlock} to keep closer relation with the problem considered.

According to \fref{BPj}, \fref{Nt-nant-eq-hspring3} and notation \fref{ospring} -- \fref{ctrac2}, we construct the expression for $\mathcal B_\omega\in\icP(\icP(\icC))^\icM$ (recall that $q$ is the identity map):
\begin{lem}
\beqnt%\label{BPjna2}
\mathcal B_\omega(\phi)=\icP\left(\bigcap_{A\in\Ntst}\bigcup_{\substack{h\in\icC\\\res hA\in\hspri{\phi}{-\omega}{A}}}\icC\res hA\right)\qquad\myll\omega\in\icOm,\ \myll\phi\in\icM.
\eeq
\end{lem}

By lemma \ref{lem-FP-unlock}, the inclusion $\mathcal F_{P_{na}}\in\UM{P_{na}}$, where the mapping $\mathcal F_{P_{na}}\in\icP(\icM)^\icM$ (see \fref{Pun}) has the form
\beq\label{FPU}
\mathcal F_{P_{na}}(\phi)\mydef\prod_{\omega\in\icOm}\icP\left(\bigcap_{A\in\Ntst}\bigcup_{\substack{h\in\icC\\\res hA\in\hspri{\phi}{-\omega}{A}}}\icC\res hA\right)\qquad\myll\phi\in\icM,
\eeq
is true.

Formally speaking the unlocking operation for the predicate $P_{na}$ is performed.
But we need some steps to apply result \fref{FPU} for solving the initial problem of constructing the greatest non-anticipating selection of the given m/f \muf.

%%%%%%%%%%%%%%%%%%%%%%%%%%%%%%%%%%%%%%%%%%%%%%%%%%%%%%%%
\subsection{Design of the greatest \Ntst-non-anticipating selection}\label{nonant-selection}

We turn to the construction of the greatest \Ntst-non-anticipating selection of $\muf$. 
So, our aim is to find $\top_{\res{P_{na}}{\icM_\muf}^{-1}(1)}$, where $\res{P_{na}}{\icM_\muf}\in\PRE{\icM_\muf}$ is the restriction of the predicate $P_{na}$ to the non-empty set $\icM_\muf\subset\icM$, $\icM_\muf\mydef\{\phi\in\icM\mid\phi\myLe\muf\}$.
By the inclusion  $\mathcal F_{P_{na}}\in\UM{P_{na}}$, we should find the greatest element among fixpoints of \fref{FPU} belonging the poset $(\icM_\muf,\myLe)$. 
One can verify that the poset $(\icM_\muf,\myLe)$ is also a complete lattice.
Hence, it is a non-empty poset.

Using lemma \ref{lem-UMP}, we construct from $\mathcal F_{P_{na}}$ an unlocking m/f  $\mathcal F_{\res{P_{na}}{\icM_\muf}}$ for the predicate \res{P_{na}}{\icM_\muf}:
\beqnt%\label{FPnaM}
\mathcal F_{\res{P_{na}}{\icM_\muf}}(\phi)=\ares{\mathcal F_{P_{na}}}{\icM_\muf}(\phi)=\icM_\muf\cap\mathcal F_{P_{na}}(\phi)=\prod_{\omega\in\icOm}\icP\left(\bigcap_{A\in\Ntst}\bigcup_{\substack{h\in\icC\\\res hA\in\hspri{\phi}{-\omega}{A}}}\muf(\omega)\res hA\right)
\eeq
for all $\phi\in\icM_\muf$.
Now we use lemma \ref{lem-f-contr} for ``narrowing'' m/f $\mathcal F_{\res{P_{na}}{\icM_\muf}}\in\UM{\res{P_{na}}{\icM_\muf}}$.
Note that the lemma is valid in our case: $(\icM_\muf,\myLe)$ is a nonempty poset. 
The m/f  $\LE_{\icM_\muf}$ (see \fref{lex}) in this case is defined by 
\beq\label{LE-in-Mf}
\LE_{\icM_\muf}(\alpha)=\prod_{\omega\in\icOm}\icP(\alpha(\omega))\qquad\alpha\in\icM_\muf.
\eeq
Following \fref{def-G} and \fref{LE-in-Mf}, we construct m/f $\bar G\in\icP(\icM_\muf)^{\icM_\muf}$: 
\beqnt%\label{lex-G}
\bar G(\phi)\mydef\mathcal F_{\res{P_{na}}{\icM_\muf}}(\phi)\cap\LE_{\icM_\muf}(\phi)
=\prod_{\omega\in\icOm}\icP\left(\bigcap_{A\in\Ntst}\bigcup_{\substack{h\in\icC\\\res hA\in\hspri{\phi}{-\omega}{A}}}\phi(\omega)\res hA\right)\qquad\myll\phi\in\icM_\muf.
\eeq
Due to inclusion $\myemp\in\icP(X)$ for any set $X$, the inequalities $\bar G(\phi)\neq\myemp$,  $\phi\in\icM_\muf$ hold. 
Consider the function $\gamma\in(\icM_\muf)^{\icM_\muf}$ defined by the rule $\gamma(\psi)\mydef\sup_{(\icM_\muf,\myLe)}\bar G(\psi)$ $\myll\psi\in\icM_\muf$: 
\begin{multline}\label{lex-g}
\gamma(\phi)=\sup\nolimits_{(\icM_\muf,\myLe)}\prod_{\omega\in\icOm}\icP\left(\bigcap_{A\in\Ntst}\bigcup_{\substack{h\in\icC\\\res hA\in\hspri{\phi}{-\omega}{A}}}\phi(\omega)\res hA\right)
\\=\prod_{\omega\in\icOm}\sup\nolimits_{(\icP(\muf(\omega)),\subset)}\icP\left(\bigcap_{A\in\Ntst}\bigcup_{\substack{h\in\icC\\\res hA\in\hspri{\phi}{-\omega}{A}}}\phi(\omega)\res hA\right)
=\prod_{\omega\in\icOm}\bigcap_{A\in\Ntst}\bigcup_{\substack{h\in\icC\\\res hA\in\hspri{\phi}{-\omega}{A}}}\phi(\omega)\res hA.
\end{multline}
for all $\phi\in\icM_\muf$.
Equalities \fref{lex-g} imply that $\gamma$ is isotone and the inclusions $\gamma(\phi)\in\bar G(\phi)$, $\myll\phi\in\icM_\muf$ are valid.
Hence, for all $\phi\in\icM_\muf$, the equality $\gamma(\phi)=\top_{\bar G(\phi)}$ is fulfilled.
Since $\bar G$ and $\gamma$ satisfy lemma \ref{lem-f-contr}, we conclude that:  
$\gamma$ is an isotone \cttv\ function and $\fix{\gamma}=P_{na}^{-1}(1)$.

The properties of the function $\gamma$ allows us to use the following theorem.

%%%%%%%%%%%%%%%%%%%%%%%%%%%%%%%%%%%%%%%%%%%%%%%%%%%%%%%
\begin{teo}[\cite{Ser_TRIMM2017e}]\label{lem-cdfix}
Let $(X,\myle)$ be a chain--complete poset, $f\in X^X$ be a \cttv\ function on $(X,\myle)$, and an ordinal $\alpha$ satisfy $\CRDP{X}\orle\alpha$.
Then $\fix{f}=\{f^\alpha(x): x\in X\}$.
\end{teo}
%%%%%%%%%%%%%%
So, for any $\alpha\in\icORD$ such that $\CRDP{\icM_\muf}\orle\alpha$, the equality
\beqnt%\label{bg-in-bG}
\res{P_{na}}{\icM_\muf}^{-1}(1)=\{\gamma^\alpha(\psi):\psi\in\icM_\muf\}
\eeq
is true.
Here we have the expression for the set of all non-anticipating selections of m/f $\muf$. 

As $\gamma$ is isotone and $\icM_\muf$ is a complete lattice, we can use the Tarski theorem \cite[Theorem 1]{Tarski1955}: the set $\fix{\gamma}=\res{P_{na}}{\icM_\muf}^{-1}(1)$ is a complete lattice in $(\icM_\muf,\myLe)$.
Hence, there is the greatest non-anticipating selection $\top_{\res{P_{na}}{\icM_\muf}^{-1}(1)}$ in the poset $\icM_\muf$. 
Due to another result of Patrick and Radhia Cousot (see \cite[Theorem 3.2]{Cousot1979}), it can be described in terms of transfinite iterations of $\gamma$ starting at \muf:
\beq\label{GNAS-iter-fix}
\top_{\res{P_{na}}{\icM_\muf}^{-1}(1)}=\top_{\fix{\gamma}}=\gamma^\alpha(\top_{\icM_\muf})=\gamma^\alpha(\muf).
\eeq
Thus, we have the desired expression for the greatest non-anticipating selection of the m/f $\muf$.

%%%%%%%%%%%%%%%%%%%%%%%%%%%%%%%%%%%%%%%%%%%%%%%%%%%%%%%%%%%%%%%%%%%%%%%%%%%%%5
\subsection{Functions $\Gamma$ and $\gamma$}

Write down expression \fref{lex-g} in the coordinate form:
\beq\label{F-eta3}
\gamma(\phi)(\omega)=\bigcap_{A\in\Ntst}\bigcup_{\substack{h\in\icC\\\res hA\in\hspri{\phi}{\omega}{A}}}\phi(\omega)\res hA\qquad\myll\omega\in\icOm\ \myll\phi\in\icM.
\eeq 
Eliminating notation \fref{ospring}, \fref{ctrac} from \fref{F-eta3}, we get the equality 
\beqnt%\label{Opnant}
\gamma(\phi)(\omega)=\Gamma(\phi)(\omega)\qquad\myll\phi\in\icM\ \myll\omega\in\icOm,
\eeq
where $\Gamma$ is given in \cite[sec. 4]{Chentsov98DUPUe}:
\beqnt%\label{Ga-Che}
\Gamma(\phi)(\omega)\mydef\{f\in\phi(\omega)\mid\myll A\in\Ntst\ \myll\omega'\in\icOm\res{\omega}{A}\myxst f'\in\phi(\omega'): \res fA=\res{f'}A\}\qquad\myll\omega\in\icOm\ \myll\phi\in\icM.
\eeq
Relation \fref{GNAS-iter-fix} generalizes the presentation \cite[theorem 6.1]{Chentsov98DUPUe}, where $\alpha=\fiord$ (the least infinite ordinal) is used.
In our case, the bigger ordinal compensates the absence of topological requirements on \icOm, \icC, and \icM.

\section{Conclusion}

The main application of the technique appears to be the existence theorems (for an equilibrium, for an equation solution).
In ordered sets, the greatest solution can be explicitly written down.

It is interesting to notice that the process (of unlocking of predicate) can also be used in the opposite direction: well known Fan's result on saddle points \cite{FanKy1953} prompted the author to look for one more fixed point theorem \cite{Ser_SUSU2016}.

\bibliographystyle{unsrt}
%%\makeatletter
%%\renewcommand{\@biblabel}[1]{#1.} % Заменяем библиографию с квадратных скобок на точку:
%%\makeatother

\bibliography{C:/Dropbox/CURRENT/ref/allbib,C:/Dropbox/CURRENT/ref/mybib}  %C:/Dropbox/CURRENT/ref/allbib

%\bibliographystyle{ifacconf} 
%\bibliography{C:/Dropbox/CURRENT/ref/allbib,C:/Dropbox/CURRENT/ref/mybib}  

\end{document}